\documentclass[10pt]{article}
\usepackage{amsmath}
\usepackage{amsfonts}
\usepackage{amssymb}
\usepackage{amsthm}
\usepackage{palatino}
\usepackage[margin=2.5cm, vmargin={1.5cm},includefoot]{geometry}

\title{On the Linear Independence of Roots}
\author{Richard Carr and Cormac O'Sullivan}
\date{September 12, 2007}
\newtheorem{theorem}{Theorem}[section]
\newtheorem{lemma}[theorem]{Lemma}
\newtheorem{prop}[theorem]{Proposition}

\begin{document}
\thispagestyle{empty}
\maketitle
\noindent

\def\H{{\mathbf{H}}}
\def\F{{\mathbb F}}
\def\C{{\mathbb C}}
\def\R{{\mathbb R}}
\def\Z{{\mathbb Z}}
\def\Q{{\mathbb Q}}
\def\N{{\mathbb N}}

\bibliographystyle{plain}

\begin{abstract}\noindent
A set of real $n$th roots that is pairwise linearly independent over the rationals must 
also be linearly independent.
We show how this result may be extended to more general fields.
\end{abstract}

\section{Introduction}
The classic Fermat equation is
\begin{equation}\label{fermatc}
x^{n} + y^{n} = z^{n}.
\end{equation}
Consider what happens when the  $n$th powers are replaced by $n$th roots
\begin{equation}\label{fermat}
x^{1/n} + y^{1/n} = z^{1/n}.
\end{equation}
Let $\N=\{1,2,\dots\}$ be the natural numbers. We seek solutions to (\ref{fermat}) with 
$x,y,z,n$ in $\N$  and $n\geqslant 2$. For simplicity we take positive
real roots and, to exclude
obvious solutions, we require that none of $x,y,z$  is a
perfect $n$th power and that $(x,y)=1$. For example, a computer search with $x,y \leqslant 
1000$ and $n \leqslant 10$ yields
\begin{equation*}
  433^{1/6} + 972^{1/6} = 42089^{1/6} + \varepsilon
\end{equation*}
with minimal error $|\varepsilon|$. It satisfies $0<|\varepsilon| < 10^{-12}$. 
Newman  shows by
elementary means in \cite{Ne} that, even
with possibly differing exponents, there are no solutions to
\begin{equation}\label{fermatb}
x^{1/m} + y^{1/n} = z^{1/r}
\end{equation}
for integers $m,n,r \geqslant 2$,  with $x,y,z$ in  $\N$,
$(x,y)=1$ and $x,y,z$ not perfect $m$th, $n$th, $r$th powers, respectively. This result
seems to have been first proven
by Obl\'ath \cite{Ob} and is also considered in
\cite{Be,XF, La, Na, To}.

\vskip 3mm
An application of our main result is to the Diophantine equation with positive rational 
exponents
 \begin{equation}\label{fermat2}
m_1 x_1^{q_1}     + m_2 x_2^{q_2} + \cdots + m_n x_n^{q_n}= 0.
\end{equation}
We are looking for solutions $(m_i,x_i,q_i)_{i=1}^n$ with
\begin{equation}\label{cond1}
m_i \in \Z,\ x_i \in \N \text{ \ and \ }0<q_i  \in \Q.
\end{equation}
 Here we
restrict to real roots, i.e. $x_i^{r/s}$ for $r,s \in \N$ means any $\alpha \in \R$ 
(possibly positive or negative)
such that $\alpha^s=x_i^r$. To avoid trivial cases we also require
\begin{equation}\label{cond2}
  m_i  \neq 0, \ x_i^{q_i} \not\in  \Z \text{ \ for each $i$  and that distinct pairs of 
}x_i
\text{s are coprime.}
\end{equation}
We will show in Proposition \ref{corol}   that solutions to
(\ref{fermat2}) satisfying (\ref{cond1}) and (\ref{cond2}) do not exist. This proposition 
follows easily
from   Theorem \ref{thm1} below. To describe it,  we first set up some notation.

\vskip 3mm
For any two fields $K \subseteq L$ define the set $\theta(K,L)$ as follows. We have $A \in 
\theta(K,L)$ if these
five conditions
are met:
\begin{enumerate}
\item $A \subseteq L$
\item $|A|  \geqslant 2$
\item For every $a\in A$ there is some  $n_a \in \N$ with $a^{n_a}\in K$. In what follows 
we always assume $n_a$ is minimal.
\item $A$ is pairwise linearly independent over $K$
\item If $\operatorname{char}(K)>0$ then $(n_a,\operatorname{char}(K))=1$ for all $a 
\in A$.
\end{enumerate}
Note that the sets $A \in \theta(K,L)$ may be infinite. What conditions on $K$ and $L$ are 
necessary so that $A \in \theta(K,L)$ is also linearly independent over $K$? For real
fields the answer is simple.

\begin{theorem} \label{thm1}
If $K\subseteq L \subseteq \R$ and $A \in \theta(K,L)$ then $A$ is linearly independent 
over $K$.
\end{theorem}
This may be generalized as follows.
\begin{theorem} \label{thm2}
If $K\subseteq L$, $A \in \theta(K,L)$ and if, for all $a \in A$, $L$ contains no $n_a$th 
root of unity except possibly $\pm 1$, then $A$ is
linearly independent over $K$.
\end{theorem}

\begin{theorem} \label{thm3}
If $K\subseteq L$, $A \in \theta(K,L)$ and if, for all $a \in A$, $K$ contains all $n_a$th 
roots of unity, then $A$ is
linearly independent over $K$.
\end{theorem}

Proposition \ref{corol} is a special case of a result  first proved by Besicovitch in 
\cite{Bes} using a type of Euclidean
algorithm for polynomials in many variables.  This proof was extended by Mordell in 
\cite{Mo} to allow the
$m_i$ and $x_i^{q_i}$ to be in more general fields. Our  Theorems \ref{thm1}, \ref{thm2}, 
\ref{thm3} provide a new
approach to these results. Their proofs are relatively short and include all cases 
considered by Besicovitch
and Mordell, see
Proposition \ref{mord}.

\vskip 3mm
A closely related question is to find the degree of the extension over $K$ you get by 
adding the roots $x_i^{q_i}$.
This was also considered in \cite{Mo} as well as in \cite{XF, Ri}. Their results are 
included in Proposition \ref{mord2}.
Siegel \cite{Si} also analyzes this question for real fields. We give a further 
application to finite fields in Proposition
\ref{ff}.

\vskip 3mm
We see from Theorems \ref{thm2}, \ref{thm3} that the roots of unity play a key role in 
these questions. The
linear dependence of roots of unity over $\Q$ is an interesting topic. For example Mann in 
\cite{Mann} proves that if
$$
m_0+m_1 \zeta^{n_1}     + m_2 \zeta^{n_2} + \cdots + m_{k-1} \zeta^{n_{k-1}}= 0
$$
for $\zeta$ a primitive $n$th root of unity, $m_i, n_i \in \Z$ and no proper subsum of the 
left side vanishing then
$$
\frac{n}{(n,n_1,n_2, \dots, n_{k-1})} \text{ \ divides \ } \prod_{p \leqslant k, p \text{ 
prime}} p.
$$
See also \cite{Co, Ev}, for example.

\section{Proof of Theorem \ref{thm1}}

We begin with the following lemma.

\begin{lemma} \label{lem1}
Let $n,\,d$ be integers with $(n,\,d)=1$ and let $\Z_n=\Z/n\Z$. Let 
$\phi:\Z_n\rightarrow\Z_n$ be given
by $\phi(x)\equiv 1+dx$ and $i(x):\Z_n\rightarrow\Z_n$ be given by $i(x)\equiv -x$.
Then the image of $\{0\}$ under iteration of $\phi,\,i$ is $\Z_n$.
\end{lemma}
\begin{proof}
Let $e$ be the inverse of $d$ so that $ed \equiv 1 \mod n$. Then $\phi$ is a bijection 
since
$\phi(e(x-1))\equiv x$.
The lemma is trivial if $d=1$ so assume $d>1$. Since $(d,n)=1$, $d$ is an element of 
the multiplicative
group $\Z_n^\ast$. Suppose $d$ has order $r$. Now
$$\phi^m(x)\equiv 1+d+d^2+\cdots+d^{m-1}+d^mx$$
and
$$(d-1)(1+d+d^2+\cdots+d^{r-1})=d^r-1\equiv 0.$$
Thus
$$
\begin{array}{rcl}
\phi^{(d-1)r}(x) & \equiv & 1+d+\cdots+d^{(d-1)r-1}+d^{(d-1)r}x\cr &
\equiv & (d-1)(1+d+d^2+\cdots+d^{r-1})+(d^r)^{d-1}x\cr & \equiv & x.
\end{array}
$$
Take the combination
$$
\begin{array}{rcl}
i(\phi(i(\phi^{(d-1)r-1})))(x) & \equiv &
-(1+d(-1-d-\cdots-d^{(d-1)r-2}-d^{(d-1)r-1}x))\cr & \equiv &
-1+d+d^2+\cdots+d^{(d-1)r-1}+d^{(d-1)r}x\cr & \equiv &
-2+1+d+d^2+\cdots+d^{(d-1)r-1}+d^{(d-1)r}x\cr & \equiv & x-2.
\end{array}
$$
So, starting from $0$ we get $-2,\,-4,\,-6,\ldots$. For $n$ odd we get all of $\Z_n$ this 
way.
If $n$ is even we get the even half of $\Z_n$, $E=\{0,\,2,\,4,\ldots,\,n-2\}$. Apply 
$\phi$ one
more time to get all of the odd elements since clearly $\phi$ is a bijection between $E$ 
and
${\Z_n}-E$ for $n$ even.
\end{proof}

With this lemma in hand, we now proceed to the proof  of Theorem \ref{thm1}.

\begin{proof}
By hypothesis, $A$ must contain a non-zero element, $a$.
Also $0$ cannot be an element of $A$. If it were, then $0 \cdot a+1 \cdot 0=0$ so that 
$\{0,\,a\}$ is linearly dependent.
Suppose now, to obtain a contradiction, that $A$ is linearly dependent over $K$. Then $A$ 
has a non-empty finite
subset which is
linearly dependent over $K$. Let $B$ be such a set of minimal cardinality. Since
 $A$ is pairwise linearly independent, the cardinality of $B$ is
not $2$. If $B$ had cardinality $1$ then it would have to be $\{0\}$ but $0$ isn't
an element of $A$ and hence of $B$. So $B$ has at least $3$ elements.
Let $I$ index $B$ so that $B=\{b_i:i\in I\}$ and $|I|\geqslant 3$. Let $K(B)$ be the field 
obtained by adjoining the elements of $B$ to $K$. This extension may not be Galois. For 
our proof to work we require a Galois extension since we will later use the fact that an 
element of a Galois extension of $K$ that is fixed by the Galois group must be in $K$.

\vskip 3mm
Let $M$ be the splitting field over $K$ of
$$
\prod\limits_{i \in I}(x^{n_{b_i}}-b_i^{n_{b_i}}).
$$
Then $M:K$ is normal and, because the characteristic of the field is $0$, separable. 
Therefore $M:K$
is Galois with Galois group $G$. It is clear that $K(B) \subseteq M$. We see that $M$ also 
contains all $n_{b_i}$th roots of unity. A short argument, left to the reader (see 
\cite[Lemma 3, p. 198]{Ros} for example), shows that $M$ contains all $n$th roots of unity 
for $n=\operatorname{lcm} \{ n_{b_i} \}_{i\in I}$ and a moment's thought reveals that $M$ 
is obtained by adjoining these roots of unity to $K(B)$.
Note too that $n$ is the minimal  number so that $b^n\in K$ for each $b\in B$.

\vskip 3mm
There are non-zero $k_i\in K$  such that $\sum_{i\in I}k_ib_i=0$.
 Let $z$ be a primitive $n$th root of unity. Clearly
$f(z)$ is again a primitive $n$th root of unity for $f\in G$, so $f(z)=z^d$ for some 
$1\leq d<n$ with $(d,n)=1$.
Also, $f(b_i)^n=f(b_i^n)=b_i^n$ so that $f(b_i)=b_iz^{t_i}$ for some $t_i$.
Let $G_d=\{f\in G:f(z)=z^d\}$ for each such $d$ so that $G$ is the disjoint union of these 
$G_d$.
For each $f\in G_d$, let $B_{t,f}=\{i\in I:f(b_i)=b_i z^t\}$ and write
$$
C_{t,f}=\sum\limits_{i\in B_{t,f}}k_ib_i$$ for $t=1, \dots ,n$.
Then the linear dependence equation is just $$\sum\limits_{t=1}^nC_{t,f}=0.$$
Apply $f$ $r$ times to get,
$$f^{r}(C_{t,f})=C_{t,f} z^{t+dt+\ldots+d^{r-1}t}.$$
Thus, for each $r$,
\begin{equation}\label{7eq}
\sum\limits_{t=1}^nC_{t,f} (z^t)^{1+d+\ldots+d^{r-1}}=0.
\end{equation}
We can also apply complex conjugation. Since each $C_{t,f}$ is real we get, for each $r$,
\begin{equation}\label{8eq}
\sum\limits_{t=1}^n C_{t,f} (z^t)^{-(1+d+\cdots+d^{r-1})}=0.
\end{equation}
With the notation of Lemma \ref{lem1} we may write (\ref{7eq}) and (\ref{8eq}) as
$$
\sum\limits_{t=1}^n C_{t,f} (z^t)^{\phi^r(0)}=0, \quad \sum\limits_{t=1}^n C_{t,f} 
(z^t)^{i(\phi^r(0))}=0
$$
respectively.
Each of these operations can be applied repeatedly.  Lemma \ref{lem1}  then implies
$$\sum\limits_{t=1}^n C_{t,f}(z^t)^k=0$$
for each $k\in\{0,\ldots,\,n-1\}$.
Let $C$ be the column vector $(C_{1,f}, \ldots ,C_{n-1,f}, C_{n,f})^T$.
Thus, we have the matrix equation $VC=0$ where $0$ is the column vector of zeros of size 
$n$
and $V$ is the Vandermonde matrix with $(i,j)$ entry
$$v_{i,j}=(z^{i-1})^j=z^{j(i-1)}.$$
This matrix is the well-known discrete Fourier transform \cite[Chapter 2]{Yip} and 
$|\det(V)|^2=n^n$. So it is invertible and $C=V^{-1} V C = V^{-1} 0 =0$.
Thus for each $t$, we have $C_{t,f}=0$. Minimality of $B$ implies that for some $t$, 
$I=B_{t,f}$.
Let $i_1,\,i_2$ be any two distinct elements of $I$.
From
$$
f\left(\frac{b_{i_1}}{b_{i_2}}\right)=\frac{f(b_{i_1})}{f(b_{i_2})}
=\frac{b_{i_1}z^t}{b_{i_2}z^t}=\frac{b_{i_1}}{b_{i_2}}
$$
we have that $b_{i_1}/b_{i_2}$ is fixed by $f$.
Now $i_1,\,i_2$ are independent of $f$ and  $d$ so that $b_{i_1}/b_{i_2}$
is in the fixed field $K$ of $G$, say $b_{i_1}/b_{i_2}=k$. Thus $1\cdot b_{i_1}-k\cdot 
b_{i_2}=0$
and $\{b_{i_1},\,b_{i_2}\}$ is linearly dependent, contradicting the assumed pairwise 
linear independence of $A$.
So, in fact, $A$ is linearly independent over $K$.
\end{proof}

\section{Generalization to arbitrary fields}
We next consider the case where the fields $K \subseteq L$ may not be real. For example, 
let $K=\Q$
and $L=\Q(A)$ for $A=\{1,\,\omega,\,\omega^2\}$, the cube
roots of unity. Then $A$ is pairwise linearly independent over $\Q$ but satisfies 
$1+\omega +\omega^2=0$.

\vskip 3mm
For another illustrative example,
consider the field $K=\Z_p(x)$  of rational functions in $x$
over $\Z_p$, the field of integers mod $p$. Let
$A=\{1,x^{1/p},(x+1)^{1/p}\}$. As we shall see, $A$ is pairwise linearly independent
over $K$  and clearly, for each $a$ in $A$, $a^p$ is in $K$. Also,
$1+x^{1/p}-(x+1)^{1/p}=0$ so $A$ is linearly dependent over $K$.
Note that $(1+x^{1/p})^p=1+x=((1+x)^{1/p})^p$ but $p$th roots are unique
as $\operatorname{char}(K)=p$.
To see the pairwise linear independence of $A$, suppose for example that
$\{1,x^{1/p}\}$ is linearly dependent. (The other cases are
similar.) So we have $f(x)^p= x g(x)^p$ for some $f(x),g(x) \in K$ with $g(x)$ not
the zero polynomial. Let $f(x)=\sum c_n x^n$. Then $f(x)^p=\sum
c_n^p x^{np}$. Similarly if $g(x)=\sum d_n x^n$ then $g(x)^p=\sum
d_n^p x^{np}$. Thus  $\sum d_n^p x^{np+1}-c_n x^{np}=0$ and $x$ is
algebraic over $\Z_p$, a contradiction.

\begin{proof}[Proof of Theorem \ref{thm2}]
We begin as before. Suppose, for a contradiction, that $A$ is linearly dependent over 
$K$.
Let $B$ be a subset of $A$ that is linearly dependent over $K$ and
minimal in cardinality with this property. Let
$B=\{b_i:i\in I\}$. As before, $B$ has at least $3$ elements. Let
$K(B)\subseteq L$ be the subfield of $L$ generated by the elements
of $B$ over $K$.  Let $M$ be the splitting field of
$$
\prod\limits_{i\in I}\left(x^{n_{b_i}}-b_i^{n_{b_i}}\right)
$$
 over $K$. Let $n=\operatorname{lcm} \{ n_{b_i} \}_{i\in I}$.
 We must have $n \geqslant 2$ since $B$ is pairwise linearly independent and
has more than one element. We see that
$M$ is also the splitting field of $x^n-1$ over $K(B)$.
As a splitting field, $M$ is normal over both $K$ and $K(B)$.
If $\operatorname{char}(K)=0$ then $M$ is also separable over both $K$ and $K(B)$. If
$\operatorname{char}(K)=p$ then note that, since $(p,n)=1$ (recall condition $(v)$ in the 
definition of
$\theta(K,L)$), each factor
$x^{n_{b_i}}-b_i^{n_{b_i}}$ is coprime to its formal derivative
$n_{b_i}x^{n_{b_i}-1}$ and so is separable.
Thus $M:K$ and $M:K(B)$ are separable and both $M:K$ and $M:K(B)$ are
Galois.
Let $z$ be a primitive
$n$th root of unity. We have the initial linear relation
\begin{equation}\label{lin}
\sum\limits_{i\in I}k_ib_i=0
\end{equation}
 where no $k_i$ is $0$. We consider
separately the cases $n=2$ and $n>2$.

\vskip 3mm
{\bf Case} $n=2$. Let $f\in \operatorname{Gal}(M:K)$. For each
$i$, $f(b_i)^2=f(b_i^2)=b_i^2$ so
$f(b_i)=c_ib_i$ with $c_i=\pm 1$.
If for some $i_1,\,i_2\in I$, we have $f(b_{i_1})=b_{i_1}$ and
$f(b_{i_2})=-b_{i_2}$ then, applying $f$ to (\ref{lin}) and adding the result to 
(\ref{lin}),
 we obtain $$
 \sum\limits_{i\in I}(1+c_i)k_i b_i=0
 $$ and since
$1+c_{i_1}=2\neq 0$ and $1+c_{i_2}=0$ we have a contradiction to the
minimality of the cardinality of $B$. So for each $f$ we have that
either $f(b_i)=b_i$ for all $i$ or that $f(b_i)=-b_i$ for all $i$.
Then for any $i_1,\,i_2\in I$, $f$ fixes $b_{i_1}/b_{i_2}$ and hence
this ratio is in the fixed field, contradicting the pairwise linear independence
of $A$ as in Theorem \ref{thm1}.

\vskip 3mm
{\bf Case} $n>2$. In
this case $z$ is not an element of $L$ by assumption. The extension $M:K(B)$ is 
cyclotomic
 and $\operatorname{Gal}(M:K(B))$ is isomorphic to $\Z_n^\ast$. Thus there is
an element $j$ of $\operatorname{Gal}(M:K(B))$ for which $j(z)=z^{-1}$. We see that $j$ 
fixes $K(B)$ and
hence $K$, so $j\in \operatorname{Gal}(M:K)$ too. Now follow the proof of  Theorem 
\ref{thm1}
but, instead of using complex conjugation, use the map $j$ to obtain
the Vandermonde matrix $V$ and demonstrate the equation $VC=0$. Again, $|\det(VV^T)|=n^n$  
and this is non-zero
since $n$ is non-zero in $K$. The rest of the proof follows as
before.
\end{proof}

\begin{proof}[Proof of Theorem \ref{thm3}]
This proof begins as in Theorems \ref{thm1} and \ref{thm2}.
Suppose, for a contradiction, that $A$ has a  subset $B$, linearly dependent over $K$, of
minimal cardinality and indexed by $I$ for $|I| \geqslant 3$. Let
$M$ be the splitting field of
$$
\prod\limits_{i\in I}\left(x^{n_{b_i}}-b_i^{n_{b_i}}\right)
$$
 over $K$. As in Theorem \ref{thm2}, $M : K$ must be Galois. The linear relation for $B$ 
is
$
\sum_{i\in I}k_ib_i=0
$
 with $k_i \neq 0$.
 Put $n=\operatorname{lcm} \{ n_{b_i} \}_{i\in I}$ and
 let $z$ be a primitive
$n$th root of unity. Then $z$ is in $K$. This requires a short argument as in the proof of 
Theorem \ref{thm1}.
Thus, for any $f \in \operatorname{Gal}(M:K)$ we have
$f(z)=z$. As in Theorem \ref{thm1}, set $B_{t,f}=\{i\in I:f(b_i)=b_i z^t\}$ and write
$
C_{t,f}=\sum_{i\in B_{t,f}}k_ib_i
$ for $t=1, \dots ,n$.
We have $\sum_{t=1}^n C_{t,f}=0$ and
Applying $f$ repeatedly shows that, for each $r \in \N$, $$\sum\limits_{t=1}^nC_{t,f} 
(z^t)^{r}=0.$$
This leads directly to the matrix equation $VC=0$ (Lemma \ref{lem1} is not required) and 
the proof continues as in
Theorem \ref{thm2}.
\end{proof}

\section{Applications}
We give some applications of Theorems \ref{thm1}, \ref{thm2} and \ref{thm3}.

\begin{prop} \label{corol}
There are no solutions to (\ref{fermat2}) satisfying (\ref{cond1}) and (\ref{cond2}).
\end{prop}
\begin{proof} Suppose that we do have a solution to (\ref{fermat2}) satisfying 
(\ref{cond1}) and (\ref{cond2}).
Take $K=\Q$ and $A=\{x_1^{q_1}, \dots , x_n^{q_n}\} \subseteq \R$. To apply
Theorem \ref{thm1} and obtain a contradiction, we need only to prove that $A \in 
\theta(\Q, \R)$ which reduces quickly to
showing that
all pairs in $A$ are linearly independent over $\Q$. If
$x_1^{r_1/s_1}$ and $x_2^{r_2/s_2}$ are linearly dependent over $\Q$, for example,
then it follows that we have $m_1,m_2 \in \Z$ with $(m_1,m_2)=1$ and
$$
m_1 x_1^{r_1/s_1} = m_2 x_2^{r_2/s_2}.
$$
Hence
$$
m_1^{s_1 s_2} x_1^{r_1 s_2} = m_2^{s_1 s_2} x_2^{r_2 s_1}.
$$
Recalling that $(x_1,x_2)=1$ we see that
$$
m_1^{s_1 s_2}=x_2^{r_2 s_1}, \ \ x_1^{r_1 s_2}=m_2^{s_1 s_2}
$$
from which we deduce that $x_1^{r_1/s_1}=\pm m_2$ and $x_2^{r_2/s_2}= \pm m_1$,
contradicting our assumption in (\ref{cond2}) that $x_i^{q_i} \not\in  \Z$.
\end{proof}

In Proposition \ref{corol} the condition in (\ref{cond2}), that the $x_i$s be pairwise 
relatively prime,
may be weakened a good deal
and  $\Q$ replaced by more general fields. This is the content of Proposition \ref{mord} 
below. For the next two results we set things up as follows.
Let $K$, $L$ be fields with $\Q \subseteq K \subseteq  L$ and let $X=\{x_1, x_2, \dots 
,x_r\}$ be a
subset of $L$ such that for every $x_i$
there is some $n_i \in \N$ (which we assume minimal) with $x_i^{n_i}\in K$.
Suppose that $X$ has the property that
$$
x_1^{e_1} x_2^{e_2} \cdots x_r^{e_r} \in K
$$
for any $(e_1, e_2, \dots, e_r ) \in \Z^r$ implies $n_i | e_i$ for all $i$ with $1 
\leqslant i \leqslant r$.
Finally, we assume that either $(i)$ $L \subseteq  \R$ or  $(ii)$ $K$ contains all $n_i$th 
roots of unity
for $1 \leqslant i \leqslant r$. Then we have the following.

\begin{prop}[Besicovitch, Mordell] \label{mord}
The $n_1 n_2 \cdots n_r$ elements $x_1^{v_1}x_2^{v_2} \cdots x_r^{v_r}$ with $0 \leqslant 
v_i < n_i$
for all $1 \leqslant i \leqslant r$ are linearly independent over $K$.
\end{prop}
\begin{proof} With Theorems \ref{thm1} and \ref{thm3} we  need only to show that
$$
A=\{ x_1^{v_1}x_2^{v_2} \cdots x_r^{v_r} \}_{0 \leqslant v_i < n_i} \in \theta(K,L).
$$
Again this reduces to proving the pairwise linear independence of elements of $A$ over 
$K$.
Take two distinct elements of $A$,
$$
a_1 = x_1^{u_1}x_2^{u_2} \cdots x_r^{u_r}, \ a_2 = x_1^{v_1}x_2^{v_2} \cdots x_r^{v_r}.
$$
If $k_1 a_1+k_2 a_2 =0$ for $k_1, k_2 \in K$ then
$$
x_1^{u_1-v_1}x_2^{u_2-v_2} \cdots x_r^{u_r-v_r} \in K
$$
and, by assumption, we have $n_i | (u_i-v_i)$ for each $i$. It follows that $a_1=a_2$ and 
this contradiction shows
that $A$ is pairwise linearly independent over $K$. Hence $A \in \theta(K,L)$
and the proof is complete.
\end{proof}

As pointed out in \cite{Mo}, Proposition \ref{mord} was also
proved by Hasse in the case where $K$ contains all $n_i$th roots of unity.

\begin{prop}[Besicovitch, Mordell] \label{mord2}
With the same notation and conditions in place we also have
$$
[K(x_1, \dots ,x_r):K]= n_1 n_2 \cdots n_r.
$$
\end{prop}
\begin{proof} With Proposition \ref{mord} we have $
[K(x_1, \dots ,x_r):K] \geqslant n_1 n_2 \cdots n_r$.  Standard results from field theory 
show the opposite inequality.
\end{proof}

Very simple proofs of the above result in the case of adjoining square roots are 
available, see \cite{Ri, Ro}.
Fried in \cite{Fr} also shows a special case of Proposition \ref{mord2} and uses it to 
give a formula for the degree of
$$
1+2^{1/2}+3^{1/3}+ \cdots +n^{1/n}
$$
over $\Q$, answering a question of Sierpi\'nski.

\vskip 3mm
Another result, based on Theorem \ref{thm1}, is the following.

\begin{prop}
Let $r_i>0$ be  roots of rationals (i.e. $r_i^{n_i} \in \Q$ with $n_i \in \N$) for $i$ 
in
a finite indexing set $I$. 
Suppose 
\begin{equation}\label{qir}
\sum
_{i \in I} q_i r_i=q \in \Q
\end{equation}
for positive rationals $q_i$. Then each $r_i$ is rational.
\end{prop}
\begin{proof}
Partition $R=\{r_i: i \in I\}$ into subsets $R_k$ for $1 \leqslant k \leqslant n$  such 
that if $\{r,r'\} \subseteq R_k$ then $\{r,r'\}$
is linearly dependent over $\Q$ and if $r,r'$ are in distinct subsets $R_k, R_l$  then
$\{r,r'\}$ is linearly independent over $\Q$.
Let $m_k$ be minimal in each $R_k$. Then $M=\{m_k : 1 \leqslant k \leqslant n \}$ is 
pairwise linearly independent over $\Q$ and hence, with Theorem \ref{thm1}, linearly 
independent over $\Q$.
But, using (\ref{qir}), the linear dependence of elements of $R_k$ and the positivity of 
$q_i, r_i$, there exist positive rationals $Q_1, \dots ,Q_n$
such that $\sum Q_k m_k=q$.
Therefore $\sum Q_k m_k+(-q)\cdot 1=0$ and $\{1\}\cup M$ is linearly dependent over $\Q$. 
Hence it is not pairwise
linearly independent over $\Q$ and some $m_k$ is rational. Without loss of generality we 
assume $k=1$, reordering if necessary.
Then
$$
\sum_{k=2}^n Q_k m_k =q-Q_1 m_1 \in \Q.
$$
If $n \geqslant 2$ then we may apply the same reasoning  to show another $m_k$ is 
rational, contradicting the pairwise linear independence of $M$.
So we must have $n=1$ and $R=R_1$  is pairwise linearly dependent over $\Q$.
Since $m_1 \in R$ is rational, it follows that all the $r_i$ are rational.
\end{proof}

\vskip 3mm
Finally, we examine how Theorems \ref{thm2} and \ref{thm3} can be used to obtain linearly 
independent sets in
 finite fields.
 Let $GF(p^u)$ denote the finite field with $p^u$ elements. If a finite field is contained 
in another,  they
necessarily
have the form $GF(p^u) \subseteq GF(p^v)$ with $u |v$. Let $m=|GF(p^u)^*|=p^u-1$ and 
$n=|GF(p^v)^*|=p^v-1$. We
put
\begin{equation}\label{l}
l=\frac nm = \frac{p^v-1}{p^u-1} =  p^{u(v/u-1)}+p^{u(v/u-2)} + \cdots +p^u +1.
\end{equation}
The next result uses the well-known fact that the multiplicative group of a finite field 
is cyclic.

\begin{prop} \label{ff}
Suppose  $GF(p^u) \subseteq GF(p^v)$ and $\phi: GF(p^v)^* \to \Z_{n}$ is an isomorphism.
Suppose also that $m>1$ and $m|l$. If the positive divisors of $m$ are $\{ d_1, d_2, 
\dots, d_r \}$ then
$$
A=\phi^{-1} \left( \left\{ \frac{d_1 l}{m}, \frac{d_2 l}{m}, \dots ,\frac{d_r l}{m} 
\right\} \right)
$$
 is linearly independent over $GF(p^u)$.
\end{prop}
\begin{proof} Note first that $\phi(GF(p^u)^*) = \langle l \rangle \subseteq \Z_{n}$. We 
verify
 that $A \in \theta(GF(p^u),GF(p^v))$. Conditions {\it (i), (ii)} are clear. If $x$ is an
element of $GF(p^v)^*$ then $x^{n_x} \in GF(p^u)$ if and only if $l | n_x \phi(x)$. It 
follows that
\begin{equation}\label{ll}
n_x=\frac{l}{(l,\phi(x))}.
\end{equation}
 Thus condition {\it (iii)} holds for all elements of $GF(p^v)$. For {\it (iv)} we can 
verify
that $x,y \in GF(p^v)^*$ are linearly independent over $GF(p^u)$ if and only if $\phi(x) 
\not\equiv \phi(y) \mod l$.
To check
{\it (v)} we need to know that $(p,n_x)=1$ for all $x \in A$. Use  (\ref{ll}) to see that 
$n_x | l$ and
 (\ref{l}) to see that $l \equiv 1 \mod p$. Thus $(p,n_x)=1$, in
fact,  for all $x \in GF(p^v)^*$. With all this $\phi^{-1}(A) \in 
\theta(GF(p^u),GF(p^v))$.

\vskip 3mm
We would like to use Theorem \ref{thm2} or \ref{thm3} to finish the proof. It may be seen 
that $GF(p^v)$ contains a $k$th
root of unity if and only if $(n,k)>1$. Since
$$
(n,n_x)=\left(n,\frac{l}{(l,\phi(x))}\right)
$$
and $l | n$ we cannot expect that $GF(p^v)$ does not contain $n_x$th roots of unity. So 
Theorem \ref{thm2} will not apply.
To use Theorem \ref{thm3} we
require $GF(p^u)$ to contain all $n_x$th roots of unity for all $x$ with $\phi(x) \in A$.
If $\zeta$ is a $k$th root of unity then $\zeta^k=1$ and $n | k \phi(\zeta)$. We see that 
all $k$ $k$th roots of unity
are in $GF(p^v)$ if and only if $k | n$ since they are
$$
\phi^{-1}\left( \left\{ 0, \frac{n}{k}, \frac{2n}{k}, \dots ,
\frac{(k-1)n}{k}\right \}\right).
$$
Clearly these are contained in $GF(p^u)$ if $l | (n/k)$ or, in other words, $k | m$. 
Therefore, with (\ref{ll}),
$GF(p^u)$ contains all $n_x$th roots of unity if
$$
\frac{l}{(l,\phi(x))} \text{ \ divides \ } m.
$$
Thus we require that
\begin{equation}\label{dlm}
\frac{l}{(l,dl/m)} \text{ \ divides \ } m
\end{equation}
for all divisors $d$ of $m$. Since $m$ also divides $l$ we have $(l,dl/m) = dl/m$. Then 
(\ref{dlm}) is easily verified and the proof is complete.
\end{proof}

For example $GF(3^2) \subseteq GF(3^{16})$ and we have $m=8$, $n=3^{16}-1$ and 
$l=5,380,840=8 \cdot w$ for $w=672,605$.
If $\phi: GF(3^{16})^* \to \Z_{n}$ is any
isomorphism then Proposition \ref{ff} implies that $\phi^{-1}(\{w,2w,4w,8w \})$ is
a subset of $GF(3^{16})$ with $4$ elements that is linearly independent over $GF(3^2)$.
Of course there exists
a set of $|GF(3^{16})/GF(3^{2})|=3^{14}$ such elements, but $A$ was found using only  
pairwise linear
independence and that
 $GF(3^2)$ contains all $8$th roots of unity.

\vskip 3mm
{\bf Acknowledgements.} We thank Gautam Chinta and Nikolaos Diamantis for their generous
help with this project.

\bibliography{linear}

\vskip .25in \noindent
Richard Carr
\newline 349, Riverside Mansions, Milk Yard, Wapping, London, E1W 3SU, UK.

\vskip .15in \noindent
Cormac O'Sullivan
\newline Department of Mathematics, Bronx Community College, Univ. Ave and W 181 St., Bronx, NY 10453, USA. \newline e-mail: cormac12@juno.com

\end{document}